\documentclass[a4paper]{amsart}
\usepackage{amsmath, amssymb, amsfonts, amscd}
\usepackage{graphicx, color}
\numberwithin{equation}{section}
\theoremstyle{plain}
 \newtheorem{thm}{Theorem}[section]
 \newtheorem{cor}[thm]{Corollary}
 \newtheorem{lem}[thm]{Lemma}
 \newtheorem{prop}[thm]{Proposition}
\theoremstyle{definition}
 \newtheorem{defn}[thm]{Definition}
 \newtheorem{exmp}[thm]{Example}

\theoremstyle{remark}
 \newtheorem{rem}[thm]{Remark}

\DeclareMathOperator{\ord}{ord}
\DeclareMathOperator{\IM}{Im}
\DeclareMathOperator{\diag}{diag}
\DeclareMathOperator{\rank}{rank}
\DeclareMathOperator{\Oidx}{Oidx}
\DeclareMathOperator{\idx}{idx}
\DeclareMathOperator{\End}{End}
\DeclareMathOperator{\spt}{spt}
\DeclareMathOperator{\id}{id}
\def\p{\partial}
\title{Katz's middle convolution and Yokoyama's extending operation}%
\author{Toshio Oshima}
\address{Graduate School of Mathematical Sciences,
University of Tokyo, 7-3-1, Komaba, Meguro-ku, Tokyo 153-8914, Japan}
\email{oshima@ms.u-tokyo.ac.jp}
\thanks{
\hspace*{12pt}Supported by Grant-in-Aid for Scientific Researches (A), 
No.\ 20244008, Japan Society of Promotion of Science\\
}
\begin{document}
\begin{abstract}
We give a concrete relation between Katz's middle convolution and Yokoyama's
extension and show the equivalence of both algorithms using these operations 
for the reduction of Fuchsian systems.
\end{abstract}

\maketitle
\section{Introduction}
Katz \cite{Kz} introduces the operations
\textsl{addition} and \textsl{middle convolution} of Fuchsian system
\begin{equation}\label{eq:Schl}
  \frac{du}{dx} = \sum_{j=1}^p \frac{A_j}{x-t_j}u
\end{equation}
of Schlesinger canonical form (SCF) on the Riemannian sphere 
and studies the rigid local systems.
It has regular singularities at $x=t_1,\dots,t_p$ and $\infty$.
Here $A_j\in M(n,\mathbb C)$ and 
$M(n,\mathbb C)$ denotes the space of $n\times n$ matrices
with entries in $\mathbb C$
and the number $n$ is called the \textit{rank} of the system.
Katz shows that any irreducible rigid system of SCF
is reduced to rank 1 system, namely a system with $n=1$, 
by a finite iteration of these operations, which implies that 
any irreducible rigid system of SCF is obtained by applying 
a finite iteration of these operations to a rank 1 system 
since these operations are invertible.  

The fact that the system is \textsl{rigid} is equal to say that it is 
free from accessory parameters but these operations are also useful
for the study of non-rigid systems.
In fact the Deligne-Simpson problem, the monodromies and integral 
representations of their solutions, their monodromy preserving 
deformations and their classification are studied by using these 
operations (cf.~\cite{DR2}, \cite{Ko}, \cite{HY}, \cite{HF}, \cite{O2}
etc.).

Dettweiler and Reiter \cite{DR} interpret these operations as those
of tuples of matrices $\mathbf A=(A_1,\dots,A_p)$ as follows.

The \textsl{addition} $M_\mu(\mathbf A)$ of $\mathbf A$ with 
$\mu=(\mu_1,\dots,\mu_p)\in\mathbb C^p$ is simply defined by
\begin{equation}
  M_\mu(\mathbf A)=M^p_\mu(\mathbf A) := (A_1 + \mu_1,\dots,A_p + \mu_p).
\end{equation}
The \textsl{convolution} $(G_1,\dots,G_p)\in M(pn,\mathbb C)^p$ of $\mathbf A$ 
with respect to $\lambda\in\mathbb C$ is define by
\begin{align}
  G_j &:= \Bigl(\delta_{\mu,j}(A_\nu+\delta_{\mu,\nu}\lambda)\Bigr)_{\substack{1\le\mu\le p\\ 1\le\nu\le p}}\qquad(j=1,\dots,p)
\allowdisplaybreaks\\
  &=
  \bordermatrix{
         &    &     &        & \underset{\smile}{j} \cr
      & \cr
    j\,{\text{\tiny$)$}} &A_1 & A_2 & \cdots & A_j+\lambda & A_{j+1} & \cdots & A_p\cr
      & \cr
  }\in M(pn,\mathbb C).\notag
\end{align}
Since the subspaces
\begin{align}
 \mathcal K&:=\bigl\{\begin{pmatrix}u_1\\ \vdots\\ u_p\end{pmatrix}\,;\,
        A_ju_j=0\quad(j=1,\dots,p)\bigr\},\allowdisplaybreaks\\
 \mathcal L_\lambda&:=\ker(G_1+\dots+G_p)
\end{align}
are $G_j$-invariant, we put $V:=\mathbb C^{pn}/\mathcal K+\mathcal L_\lambda$
and define $\bar G_j\in\End(V)\simeq M(\dim V,\mathbb C)$, which are
the linear maps induced by $G_j$, respectively.
Then the \textsl{middle convolution} $mc_\lambda(\mathbf A)$ of $\mathbf A$
equals $(\bar G_1,\dots,\bar G_p)$.

For $\mathbf A=(A_1,\dots,A_p)$, 
$\mathbf B=(B_1,\dots,B_p)\in M(n,\mathbb C)^p$
we write $\mathbf A\sim\mathbf B$ if there exists $g\in GL(n,\mathbb C)$
such that $B_j=gA_jg^{-1}$ for $j=1,\dots,p$ and we will sometimes
identify $\mathbf A$ with $\mathbf B$ if $\mathbf A\sim\mathbf B$.
The corresponding systems \eqref{eq:Schl} will be also identified.

Yokoyama \cite{Yo} introduces an extension and a restriction of the 
Fuchsian system
\begin{equation}\label{eq:Okubo}
  (xI_n-T)\frac{du}{dx}=Au
\end{equation}
of Okubo normal form (ONF) with $A$, $T\in M(n,\mathbb C)$ when $T$ is a 
diagonal matrix and $A$ satisfies a certain condition.

Suppose
\begin{equation}\label{eq:T}
  T = \begin{pmatrix}
       t_1I_{n_1}\\
       &\ddots\\
       &&t_pI_{n_p}
      \end{pmatrix},
\end{equation}
where $n=n_1+\dots+n_p$ is a partition of $n$ and
$t_i\ne t_j$ if $i\ne j$.  Put
\begin{equation}\label{eq:O2S}
  A=\begin{pmatrix}
    A_{11}&\cdots&A_{1p}\\
    \vdots&\cdots&\vdots\\
    A_{p1}&\cdots&A_{pp}
    \end{pmatrix}
\end{equation}
according to the partition, namely $A_{ij}\in M(n_i,n_j;\mathbb C)$
which is the space of $n_i\times n_j$ matrices with entries 
in $\mathbb C$.
Here we note that the system \eqref{eq:Okubo} of ONF
is equal to the system \eqref{eq:Schl}
of SCF by putting\\[-.9cm]
\begin{align}
  A_j:= \bordermatrix{
      &\cr
      & \cr
    j\,{\text{\tiny$)$}} &A_{j1} & A_{j2} & \cdots & A_{jp}\cr
      & \cr
  }\in M(n,\mathbb C).\label{eq:OkSh}
\end{align}
Conversely we have the following lemma.
\begin{lem}\label{lem:Okubo}
Suppose $(A_1,\dots,A_p)\in M(n,\mathbb C)^p$ satisfies
\begin{equation}\label{C:Okubo}
  \begin{cases}
   \rank A_1+\cdots+\rank A_p = n,\\
   \IM A_1+\cdots+\IM A_p =\mathbb C^n.
  \end{cases}
\end{equation}
Then there exists $g\in GL(n,\mathbb C)$ such that
the $\nu$-th row of $g^{-1}A_jg$ is identically zero
if $\nu\le \rank A_1+\cdots+\rank A_{j-1}$ or 
$\nu>\rank A_1+\cdots+\rank A_{j}$. 
Hence the system of SCF is equivalent to a system of ONF
if \eqref{C:Okubo} holds.
\end{lem}
\begin{proof}
The assumption implies that there exists a basis $\{v_1,\dots,v_n\}$
of $\mathbb C^n$ such that
\[
  \IM A_j=
  \sum_{\rank A_1+\cdots+\rank A_{j-1}<\nu\le \rank A_1+\cdots+\rank A_j}
  \mathbb Cv_j.
\]
Then the expression of $A_j$ under this basis has the required
property, namely, we may put $g=(v_1,\dots,v_n)\in GL(n,\mathbb C)$.
\end{proof}
\begin{rem}\label{rem:C:Okubo}
If a system \eqref{eq:Okubo} of ONF is linearly irreducible
(cf.~Definition~\ref{def:linirred}), it satisfies 
\eqref{C:Okubo} with \eqref{eq:O2S} and \eqref{eq:OkSh}.
\end{rem}

Yokoyama \cite{Yo} defines \textsl{extentions} 
$(\hat T,\hat A)=E_\epsilon(T,A)$ for $\epsilon=0,1$ and $2$
with respect to two distinct complex numbers $\rho_1$ and $\rho_2$ 
when $T$ and $A_{ii}$ ($i=1,\dots,p$) are diagonalizable. 
Here $\epsilon$ is the number of the elements of $\{\rho_1,\rho_2\}$ which
are not the eigenvalues of $A$.

Let
\begin{align}
 A_{ii}\sim \begin{pmatrix}
         \lambda_{i,1}I_{\ell_{i,1}}\\
            &\ddots\\
          &&\lambda_{i,r_i}I_{\ell_{i,r_i}}
        \end{pmatrix}\label{eq:Aii}
\end{align}
with $\lambda_{i,j}\ne\lambda_{i,k}$ ($j\ne k$) and
$n_i=\ell_{i,1}+\cdots+ \ell_{i,r_i}$ and
fix a matrix $P\in GL(n,\mathbb C)$ so that
\begin{align}
  A':=\begin{pmatrix}
      \mu_1 I_{m_1}\\
      &\ddots\\
      &&\mu_q I_{m_q}
     \end{pmatrix}=P^{-1}AP\sim A,\label{eq:A}
\end{align}
where $n=m_1+\cdots+m_q$ and $\mu_i\ne\mu_j$ ($i\ne j$).
Then $E_2(T,A)=(\hat T,\hat A)$ with
\begin{align}
 \hat T&:=\begin{pmatrix}
          T\\& t_{p+1}I_n
         \end{pmatrix},\\
 \hat A&:=\begin{pmatrix}
          A & P\\
          -(A'-\rho_1I_n)(A'-\rho_2I_n)P^{-1} & (\rho_1+\rho_2)I_n - A'
         \end{pmatrix}.
\end{align}
When $\rho_1$ or $\rho_2$ is an eigenvalue of $A$, there exists
a subspace invariant by $\hat T$ and $\hat A$ and the 
\textsl{extending operations}
$E_1$ and $E_0$ of $(T,A)$ are defined as follows.
Putting
\begin{align}
 V_k:=\bigl\{\begin{pmatrix}
        u\\ v_1 \\ v_2
      \end{pmatrix}\,;\,u\in \mathbb C^n,\ v_1=0\in\mathbb C^k
      \text{ and } 
      v_2\in\mathbb C^{n-k}\bigr\}.
\end{align}
we have
\begin{align}
E_1(T,A)&:=(\hat T|_{V_{m_1}},\hat A|_{V_{m_1}})
&&\text{when \ }\rho_1=\mu_1,\\
E_0(T,A)&:=(\hat T|_{V_{m_1+m_2}},\hat A|_{V_{m_1+m_2}})
&&\text{when $\rho_1=\mu_1$ and $\rho_2=\mu_2$}.
\end{align}
\textsl{Restrictions} are defined as inverse operations of these
extensions.
It is proved by \cite{Yo} that any irreducible rigid system of
ONF with generic spectral parameters 
$\lambda_{i,j}$ and $\mu_k$ is reduced to a rank 1 system 
by a finite iteration of the extensions and restrictions
and it gives the monodromy of the system.

In this note we clarify the direct relation between Yokoyama's operations
and Katz's operations and then relax the assumption to define 
Yokoyama's operations (cf.~Theorem~\ref{thm:ext} and Theorem~\ref{thm:rest}).
In particular we don't assume that the local monodromies of the system 
are semisimple (cf.~Theorem~\ref{thm:Oreg}).
Moreover we show in Theorem~\ref{thm:eq} that the both operations on Fuchsian 
systems are equivalent in a natural sense.
Hence the property of Katz's operation can be transferred to that of
Yokoyama's operations and vice versa.
For example, it is proved by \cite{HF} that the middle convolution 
preserves the deformation equation and therefore so do Yokoyama's operations.

The author would like to express his sincere gratitude to
Y.~Haraoka and T.\ Yokoyama for valuable discussions.
\section{Katz's middle convolution}
For a partition $\mathbf m=(m_1,\dots,m_N)$ of $n$ with $n=m_1+\cdots+m_N$
and $\lambda=(\lambda_1,\dots,\lambda_N)\in\mathbb C^N$
we define a matrix $L(\mathbf m;\lambda)\in M(n,\mathbb C)$ as  
a representative of a conjugacy class, which is 
introduced and effectively used by \cite{Os}
(cf.~\cite[\S3]{O2}):

If $m_1\ge m_2\ge\cdots\ge0$, then
\begin{equation}\label{eq:Ln}
\begin{split}
 L(\mathbf m;\mathbf \lambda) 
  &:= \Bigl(A_{ij}\Bigr)_{\substack{1\le i\le N\\1\le j\le N}},\quad
 A_{ij}\in M(m_i,m_j,\mathbb C),\\
 A_{ij} &= \begin{cases}
          \lambda_i I_{m_i}&(i=j)\\
          I_{m_i,m_j}:=
          \Bigl(\delta_{\mu\nu}\Bigr)
          _{\substack{1\le \mu\le m_i\\1\le \nu\le m_j}}
	=
          \begin{pmatrix}
          I_{m_j} \\ 0
          \end{pmatrix}&(i=j-1)\\
          0            &(i\ne j,\ j-1)
          \end{cases}.
\end{split}
\end{equation}
For example
\begin{equation*}
 L(2,1,1;\lambda_1,\lambda_2,\lambda_3)=
 \begin{pmatrix}
 \lambda_1 & 0       & 1& 0\\
 0         &\lambda_1& 0& 0\\
 0         & 0       &\lambda_2&1\\
 0         & 0       & 0       &\lambda_3\\
 \end{pmatrix}.
\end{equation*}
Denoting $
  Z_{M(n,\mathbb C)}(A):=
  \{X\in M(n,\mathbb C)\,;\,AX=XA\}$, we have
\begin{align}
  \dim\ker \prod_{j=1}^k\bigl(L(\mathbf m;\lambda)-\lambda_j\bigr)
    &=m_1+\cdots+m_k
 \quad(k=1,\dots,N),\\
  \dim Z_{M(n,\mathbb C)}\bigl(L(\mathbf m;\lambda)\bigr)&=
  m_1^2+\cdots+m_N^2.\label{eq:ZL}
\end{align}
In general we fix a permutation $\sigma$ of indices $1,\dots,N$
so that $m_{\sigma(1)}\ge m_{\sigma(2)}\ge\cdots$ and define
$L(\mathbf m;\lambda)=L(m_{\sigma(1)},\dots,m_{\sigma(N)};
\lambda_{\sigma(1)},\dots,\lambda_{\sigma(N)})$.

Let $\mathbf A=(A_1,\dots,A_p)\in M(n,\mathbb C)^p$.
Put
\begin{equation}\label{eq:A0}
  A_0=-(A_1+\cdots+A_p).
\end{equation}
Then Katz \cite{Kz} defines
\begin{align}
  \idx\mathbf A:=\sum_{j=0}^p\dim Z_{M(n,\mathbb C)}(A_j)
  - (p-1)n^2,
\end{align}
which is called the \textsl{index of rigidity}.

If $\mathbf A$ is irreducible, $\idx\mathbf A\le 2$.
Moreover an irreducible $\mathbf A$ is rigid if and only if 
$\idx\mathbf A=2$, which is proved by \cite[\S1.1.1]{Kz}.
Here $\mathbf A$ is called \textsl{irreducible} if
any subspace $V$ of $\mathbb C^n$ satisfying
$A_j V\subset V$ for $j=1,\dots,p$ is $\{0\}$ or $\mathbb C^n$.

Using the representatives $L(\mathbf m;\lambda)$ of conjugacy classes of
matrices, we can easily describes the property of the 
middle convolution.
\begin{defn}
For $\mathbf A\in M(n,\mathbb C)^p$ we 
choose a tuple of $p+1$ partitions
$\mathbf m=(\mathbf m_0,\dots,\mathbf m_p)$ and 
$\lambda_{j,\nu}\in\mathbb C$ so that
\begin{equation}\label{eq:DS}
 A_j\sim L\bigl(\mathbf m_j;\lambda_j)
 \text{ with }\mathbf m_j:=(m_{j,1},\dots,m_{j,n_j})\text{ and }
 \lambda_j:=(\lambda_{j,1},\dots,\lambda_{j,n_j})
\end{equation}
for $j=0,\dots,p$.
Here $A_0$ is determined by \eqref{eq:A0}.
We define the \textsl{Riemann scheme} of the corresponding system 
\eqref{eq:Schl} of SCF by
\begin{equation}\label{eq:RS}
 \begin{Bmatrix}
    x=\infty& x=t_1 & \cdots & x=t_p \\
   [\lambda_{0,1}]_{(m_{0,1})}& [\lambda_{1,1}]_{(m_{1,1})}&\cdots&
     [\lambda_{p,1}]_{(m_{p,1})}\\
   \vdots&\vdots&   &\vdots\\
     [\lambda_{0,n_0}]_{(m_{0,n_0})}& [\lambda_{1,n_1}]_{(m_{1,n_1})}&\cdots&
     [\lambda_{p,n_p}]_{(m_{p,n_p})}
 \end{Bmatrix}.
\end{equation}
Here $A_j$ is called the \textsl{residue matrix} of the system at 
$x=t_j$ ($j=1,\dots,p$) and $A_0$ is the residue matrix of the system
at $x=\infty$.
We also call \eqref{eq:RS} the Riemann scheme of $\mathbf A$.
We will allow that some $m_{j,\nu}$ are $0$.
\end{defn}
\begin{thm}[\cite{DR}, \cite{DR2}]\label{thm:mid}
Let $\mathbf A=(A_1,\dots,A_p)\in M(n,\mathbb C)^p$.
Assume the following conditions:
\begin{align}
  \bigcap_{\substack{1\le\nu\le p\\ \nu\ne i}}\ker A_\nu
  \cap \ker(A_i+\tau)&=\{0\}&(i=1,\dots,p,\ \forall\tau\in\mathbb C),
  \label{eq:star}\\
  \sum_{\substack{1\le\nu\le p\\ \nu\ne i}}\IM A_\nu
  +\IM(A_i+\tau)&=\mathbb C^n&(i=1,\dots,p,\ \forall\tau\in\mathbb C)
  \label{eq:starstar}.
\end{align}
Then $\bar{\mathbf G}=(\bar G_1,\ldots,\bar G_p)
:=mc_{\lambda}(\mathbf A)$ satisfies \eqref{eq:star} and
\eqref{eq:starstar} and
\begin{equation}
 \idx \bar{\mathbf G} =\idx\mathbf A.
\end{equation}
If $\mathbf A$ is irreducible, so is $\bar{\mathbf G}$.  
If $\mathbf A\sim\mathbf B$, 
then $mc_\lambda(\mathbf A)\sim mc_\lambda(\mathbf B)$.
Moreover we have
\begin{align}
 mc_0(\mathbf A)&\sim \mathbf A,\\
 mc_{\lambda'}\circ mc_\lambda(\mathbf A)&\sim mc_{\lambda'+\lambda}
(\mathbf A).
\end{align}

Let \eqref{eq:RS} be the Riemann scheme of $\mathbf A$.
We may assume
\begin{equation}
 \begin{cases}
   \lambda_{0,1} = \lambda,\\
   \lambda_{i,0} = 0 &(i=1,\dots,p),\\
   \lambda_{j,\nu} = \lambda_{j,0}\ \Rightarrow\ m_{j,\nu}\le m_{j,0}
                     &(\nu=1,\dots,n_j,\ j=0,\dots,p).
 \end{cases}
\end{equation}
Note that $m_{j,0}$ may be $0$.
Then the Riemann scheme 
\begin{align}
 \begin{Bmatrix}
    x=\infty& x=t_1 & \cdots & x=t_p \\
   [\lambda]_{(m_{0,1})}& [0]_{(m_{1,1})}&\cdots&
     [0]_{(m_{p,1})}\\
   [\lambda_{0,2}]_{(m_{0,2})}& [\lambda_{1,2}]_{(m_{1,2})}
     &\cdots&[\lambda_{p,2}]_{(m_{p,2})}\\
   \vdots&\vdots&   &\vdots\\
     [\lambda_{0,n_0}]_{(m_{0,n_0})}&
     [\lambda_{1,n_1}]_{(m_{1,n_1})}
     &\cdots&[\lambda_{p,n_p}]_{(m_{p,n_p})}
 \end{Bmatrix}
\end{align}
of $\mathbf A$ is transformed into the Riemann scheme
\begin{align}
 \begin{Bmatrix}
    x=\infty& x=t_1 & \cdots & x=t_p \\
   [-\lambda]_{(m_{0,1}-d)}& [0]_{(m_{1,1}-d)}&\cdots&
     [0]_{(m_{p,1}-d)}\\
   [\lambda_{0,2}-\lambda]_{(m_{0,2})}& [\lambda_{1,2}+\lambda]_{(m_{1,2})}
     &\cdots&[\lambda_{p,2}+\lambda]_{(m_{p,2})}\\
   \vdots&\vdots&   &\vdots\\
     [\lambda_{0,n_0}-\lambda]_{(m_{0,n_0})}&
     [\lambda_{1,n_1}+\lambda]_{(m_{1,n_1})}
     &\cdots&[\lambda_{p,n_p}+\lambda]_{(m_{p,n_p})}
 \end{Bmatrix}
\end{align}
of $mc_\lambda(\mathbf A)$ with
\begin{equation}
d=m_{0,1}+\cdots+m_{p,1}-(p-1)n.
\end{equation}
\end{thm}


\begin{rem} \label{rem:irred}
If $\mathbf A$ is irreducible, then \eqref{eq:star} and \eqref{eq:starstar}
are valid.
\end{rem}

Suppose $\lambda\ne0$.
Since
\begin{multline}
\begin{pmatrix}
   A_1\\
    & \ddots\\
    && A_p
 \end{pmatrix}
 \begin{pmatrix}
   A_1+\lambda&\cdots&A_p\\
   \vdots & \cdots&\vdots\\
   A_1& \cdots & A_p+\lambda
 \end{pmatrix}\\
 =
  \begin{pmatrix}
   A_1+\lambda&\cdots&A_1\\
   \vdots & \cdots&\vdots\\
   A_p& \cdots & A_p+\lambda
 \end{pmatrix}
 \begin{pmatrix}
   A_1\\
    & \ddots\\
    && A_p
 \end{pmatrix}\label{eq:midIM}
\end{multline}
and the linear map $A_j$ induces the isomorphism
$\mathbb C^n/\ker A_j\simeq\IM A_j$,
we put
\begin{align}
\tilde A&:=
  \begin{pmatrix}
   A_1\\
    & \ddots\\
    && A_p
 \end{pmatrix},\\
G'_j&:=\bordermatrix{
         &    &     &        & \underset{\smile}{j} \cr
      & \cr
    j\,{\text{\tiny$)$}} &A_j& A_j & \cdots & A_j+\lambda & A_j & \cdots & A_j\cr
      & \cr
  }\in M(pn,\mathbb C)
\end{align}
for $j=1,\dots,p$ and define
\begin{equation}
 \begin{split}
  G'_0&:=-(G'_1+\cdots+G'_p),\\
  \bar G'_j&:=
      G'_j\!\bigm|_{\IM\tilde A/\ker G'_0}
  \qquad(j=0,\dots,p),\\
 \end{split}
\end{equation}
\begin{lem}\label{lem:mcbyIM}
Suppose $\lambda\ne0$ and put $G_0=-(G_1+\cdots+G_p)$.
Then under the above notation
\begin{align}
 \tilde A G_j&=G'_j \tilde A&&(j=0,\dots,p),\label{eq:1-1}\\
 G'_j\IM\tilde A&\subset\IM\tilde A,\quad
 \ker G'_0=\bigcap_{j=1}^p\ker G'_j &&(j=0,\dots,p),\label{eq:1-2}\\
 \tilde A(\mathcal K+\mathcal L_\lambda)&=\ker G'_0
 \subset\IM\tilde A,\label{eq:1-3}
\end{align}
and therefore $\tilde A\in\End(\mathbb C^{pn})$ induces
the isomorphism
\begin{equation}
\begin{aligned}
 (\bar G_1,\dots,\bar G_p)&=mc_\lambda(A_1,\dots,A_p)
 &&\in\Bigl(\End\bigl(\mathbb C^{pn}\!/\mathcal K+\mathcal L_\lambda\bigr)\Bigr)^p\\
  &\sim (\bar G'_1,\dots,\bar G'_p)
  &&\in\Bigl(\End\bigl(\IM\tilde A/\ker G'_0\bigr)\Bigr)^p.
\end{aligned}
\end{equation}

In particular if $-\lambda$ is not the eigenvalue of $A_1+\cdots+A_p$,
the middle convolution $mc_\lambda(\mathbf A)$ transforms the system 
\eqref{eq:Schl} of SCF to the system of ONF
\begin{equation}\label{eq:2Okubo}
 \left(xI_{n'_1+\cdots+n'_p} - \begin{pmatrix}
             t_1I_{n'_1}\\
             &\ddots\\
             &&t_pI_{n'_p}
           \end{pmatrix}\right)\frac{du}{dx}
 =\bigl(-G'_0\!\bigm|_{\IM A_1\oplus\cdots\oplus\IM A_p}\bigr)u
\end{equation}
with $n'_j=\dim\IM A_j$.
\end{lem}
\begin{proof}
Note that $\tilde AG_0=G'_0\tilde A$, which corresponds to \eqref{eq:midIM},
and moreover that \eqref{eq:1-1} and \eqref{eq:1-2} are also clear.

Since $\mathcal K=\ker\tilde A$ and $\mathcal L_\lambda=\ker G_0$,
$G'_0\tilde A(\mathcal K+\mathcal L_\lambda)=G'_0\tilde A\ker G_0
=\tilde AG_0\ker G_0=0$ and therefore
$\tilde A(\mathcal K+\mathcal L_\lambda)\subset\ker G'_0$.
Let $u\in\ker G'_0$.  Putting
\[
 u=\begin{pmatrix}u_1\\ \vdots\\ u_p\end{pmatrix},\quad
 u_j\in\mathbb C^n,\quad
 v:=u_1+\dots+u_p\text{ \ and \ }
 \tilde v:=\begin{pmatrix}v\\ \vdots\\ v\end{pmatrix},
\]
we have $\lambda u_j=-A_jv$
and therefore $\lambda v+(A_1+\cdots+A_p)v=0$.
Hence
$\tilde v\in\ker G_0$ and $u=-\lambda^{-1}\tilde A\tilde v\in\tilde A\ker G_0$.
Thus we have \eqref{eq:1-3}.
\end{proof}
\section{Yokoyama's extending operation}
First we examine the conditions \eqref{eq:star} and \eqref{eq:starstar}
for the Fuchsian system \eqref{eq:Okubo} of ONF with \eqref{eq:T}.

For a partition $n=k_1+\cdots+k_q$ and $C_j\in M(k_j,\mathbb C)$
we denote
\begin{align*}
  \diag(C_1,\dots,C_q)&:=
  \begin{pmatrix}
   C_1\\  &\ddots\\ && C_p
  \end{pmatrix}\in M(n,\mathbb C),\\
  O_{k_j} &:= 0\in M(k_j,\mathbb C).
\end{align*}
Then $A_j$ given by \eqref{eq:OkSh} equals 
$\diag(O_{n_1+\cdots+n_{j-1}},I_{n_j},O_{n_{j+1}+\cdots+n_p})A$.
\begin{lem}\label{lem:Ocond}
The pair of conditions \eqref{eq:star} and  \eqref{eq:starstar}
for $A_j$ given by \eqref{eq:OkSh}
is equivalent to the pair of conditions
\begin{gather}
 \rank A=n\label{eq:rank_n}
 \intertext{and}
 \begin{cases}
  \rank \bigl((A+\tau)
  \diag(O_{n_1+\cdots+n_{i-1}},I_{n_i},O_{n_{i+1}+\cdots+n_p})\bigr)
  = n_i,\\ 
 \rank \bigl(\diag(O_{n_1+\cdots+n_{i-1}},I_{n_i},O_{n_{i+1}+\cdots+n_p})
  (A+\tau)\bigr)
 = n_i
 \end{cases}\label{eq:Ostar}\\
 \text{\hspace{1cm}for any $\tau\in\mathbb C$ and $j=i,\dots,p$.}\notag
\end{gather}
\end{lem}
\begin{proof}
Note that the condition \eqref{eq:star} with $\tau=0$ equals \eqref{eq:rank_n},
which implies \eqref{eq:Ostar} with $\tau=0$.

Suppose $\tau\ne0$ and \eqref{eq:rank_n}. Put
$\mathbf u=\begin{pmatrix}u_1\\ \vdots\\ u_p\end{pmatrix}$
with $u_\nu\in\mathbb C^{n_\nu}$.
Then $\sum_{\nu\ne i}\IM A_\nu
=\{\mathbf u\in\mathbb C^n\,;\,u_i=0\}$ and therefore
the condition \eqref{eq:starstar} is equivalent to
the second condition of \eqref{eq:Ostar}.
Since $\ker(A_i+\tau)=\{\mathbf u\in\mathbb C^n\,;\,
(A_{ii}+\tau)u_i=0\text{ and }u_\nu=0\quad(\nu\ne i)\}$,
the condition \eqref{eq:star} is equivalent to the condition
$\{u_i\in\mathbb C^{n_i}\,;\,(A_{ii}+\tau)u_i=0\text{ and }
A_{\nu,i}u_i=0\quad(\nu\ne i)\}=\{0\}$, which is equivalent to the 
second condition of \eqref{eq:Ostar}.
%
%
%
\end{proof}
\begin{defn}\label{def:linirred}
The system \eqref{eq:Schl} of SCF is called \textsl{linearly irreducible}
if $A_j$ have no non-trivial common invariant subspace of $\mathbb C^n$,
namely, $\mathbf A=(A_1,\dots,A_p)$ is irreducible.
Then
\begin{equation}
\text{irreducible\ }\Rightarrow\text{\ linearly irreducible\ }
\Rightarrow\text{\ \eqref{eq:rank_n} and \eqref{eq:Ostar}}.
\end{equation}
\end{defn}
\begin{rem}
The conditions \eqref{eq:rank_n} and \eqref{eq:Ostar} are valid 
if the system \eqref{eq:Okubo} of ONF is irreducible as a differential equation
or linearly irreducible.
\end{rem}

Assume \eqref{eq:rank_n} and \eqref{eq:Ostar} for the system \eqref{eq:Okubo}
of ONF.
Put $\lambda=-\rho_1\ne0$ and
apply Lemma~\ref{lem:mcbyIM} to $\mathbf A=(A_1,\dots,A_p)$ given by
\eqref{eq:OkSh}.
Then \eqref{eq:rank_n} assures $\IM A_j\simeq \mathbb C^{n_j}$.
Under the notation in the proof Lemma~\ref{lem:Ocond} 
the projection defined by
\[
 \iota_j:\mathbb C^n\ni
 \mathbf u=\begin{pmatrix}u_1\\ \vdots\\ u_p\end{pmatrix}\mapsto
 u_j\in \mathbb C^{n_j}
\]
gives this isomorphism and hence we have the isomorphism 
$\iota:\IM\tilde A\simeq\mathbb C^{n_1+\cdots+n_p}
=\mathbb C^n$.
Under the identification of this isomorphism $\iota$ we have
\begin{align*}
 G'_j|_{\IM\tilde A}&\simeq 
 G''_j:=\bordermatrix{
         &    &     &        & \underset{\smile}{j} \cr
      & \cr
    j\,{\text{\tiny$)$}} &A_{j1}& A_{j2} & \cdots & A_{jj}-\rho_1
       & A_{jj+1} & \cdots & A_{jp}\cr
      & \cr
  }\in M(n,\mathbb C),\\
  G''_1&+\cdots+G''_p=A-\rho_1,\\
 G'_j&|_{\IM\tilde A/\ker G'_0}\simeq
  \bar G''_j:=G''_j|_{\mathbb C^n\!/\ker(A-\rho_1)}
\end{align*}
for $j=1,\dots,p$ and
\begin{multline}
mc_{-\rho_1}(A_1,\dots,A_p)\in
 \Bigl(\End\bigl(\mathbb C^{pn}\!/\mathcal K+\mathcal L_\lambda\bigr)\Bigr)^p\\
  \sim (\bar G''_1,\dots,\bar G''_p)
 \in\Bigl(\End\bigl(\mathbb C^n\!/\ker G'_0\bigr)\Bigr)^p.
\end{multline}
In particular we have
\begin{cor}\label{cor:Euler}
Suppose the system \eqref{eq:Okubo} of ONF satisfies \eqref{eq:rank_n} and 
\eqref{eq:Ostar}.
If $-\lambda$ is not the eigenvalue of $A$, then
the middle convolution $mc_{\lambda}(A_1,\dots,A_p)$ corresponds to the
transformation $A\mapsto A+\lambda$ of the system \eqref{eq:Okubo}.
\end{cor}
\begin{defn}
We denote this operation of the system of ONF by $E_\lambda$ and call it
a \textsl{generic Euler transformation}, which is defined if $-\lambda$ 
is not the eigenvalue of $A$.  
Note that $E_\lambda\circ E_{\lambda'}=E_{\lambda+\lambda'}$.
\end{defn}
The transformation $A\mapsto A+\lambda$ of \eqref{eq:Okubo}
corresponds the Riemann-Liouville integral
\begin{equation}
  I^\lambda_t u(x):=\frac1{\Gamma(\lambda)}\int_t^x(x-s)^{\lambda-1}u(s)ds
\end{equation}
of the solution $u(x)$ of the system (cf.~\cite[Chapter 5]{Kh}).
Here $t\in\{t_1,\dots,t_p,\infty\}$.
\begin{defn}\label{def:T}
Define the linear maps
\begin{equation*}
 \begin{matrix}
   T_{(j,\infty)}:& M(n,\mathbb C)^p &\to& M(n,\mathbb C)^p\\
     &\text{\rotatebox{90}{$\in$}}&&\text{\rotatebox{90}{$\in$}}\\
     & \bigl(B_1,\dots,B_p\bigr) &\mapsto&
       \bigl(B_1,\dots,B_{j-1},-(B_1+\dots+B_p),B_{j+1},\dots,B_p\bigr)
 \end{matrix}
\end{equation*}
for $j=1,\dots,p$ and
\begin{equation*}
 \begin{matrix}
  T_\sigma:& M(n,\mathbb C)^p &\to& M(n,\mathbb C)^p\\
      &\text{\rotatebox{90}{$\in$}}&&\text{\rotatebox{90}{$\in$}}\\
      & \bigl(B_1,\dots,B_p\bigr) &\mapsto&
        \bigl(B_{\sigma(1)},\dots,B_{\sigma(p)}\bigr).
 \end{matrix}
\end{equation*}
Here $\sigma$ is a permutation of the indices $1,\dots,p$.
Under the natural identification
\begin{equation}
 M(n,\mathbb C)^p\simeq \{(B_1,\dots,B_{p+1})\in M(n,\mathbb C)^{p+1}\,;\,
 B_{p+1}=0\}\subset M(n,\mathbb C)^{p+1}
\end{equation}
we have $T_{(p+1,\infty)}(B_1,\dots,B_p)=(B_1,\dots,B_p,-(B_1+\dots+B_p))$.
\end{defn}
\begin{rem}
{\rm i) \ }
Let $\mathbf B\in M(n,\mathbb C)^p$.
Then $T_{(p+1,\infty)} \mathbf B$ is irreducible if and only if $\mathbf B$
is irreducible.
Similarly $T_{(p+1,\infty)} \mathbf B$ satisfies \eqref{eq:star} 
and \eqref{eq:starstar} if and only if so does $\mathbf B$.

{\rm ii) \ } The map $T_{(p+1,\infty)}$ corresponds to the 
transformation of the Fuchsian system of SCF induced from 
the automorphism of the Riemannian
sphere defined by $x\mapsto \frac{t_{p+1}x-c}{x-t_{p+1}}$.
Here $c\in\mathbb C$, $c\ne t_{p+1}^2$ and $t_{p+1}\ne t_j$ for $j=1,\dots,p$.

{\rm iii) } 
The middle convolution $mc_\lambda$ clearly commutes with $T_\sigma$,
namely,
\begin{align}
  mc_\lambda\circ T_\sigma = T_\sigma\circ mc_\lambda.
\end{align}
\end{rem}

Fix $\rho_2\ne0$ and examine 
$mc_{\rho_1}\circ
 M_{(0,\dots,0,\rho_2-\rho_1)}\circ
 T_{(p+1,\infty)}\circ
 mc_{-\rho_1}(A_1,\dots,A_p)$.

Since $G''_j\ker(A-\rho_1)=0$, it follows from \eqref{eq:Ostar} that 
$\IM\bar G''_j=\IM G''_j/G''_j\ker(A-\rho_1)\simeq\mathbb C^{n_j}$.
Note that
\begin{align*}
 M_{(0,\dots,0,\rho_2-\rho_1)}\circ
 T_{(p+1,\infty)}(\bar G''_1,\dots,\bar G''_p)
 &=(\bar G''_1,\dots,\bar G''_p,-\bar G''_1-\cdots-\bar G''_p+\rho_2-\rho_1)\\
 &=\bigl(\bar G''_1,\dots,\bar G''_p,(-A+\rho_2)\!
  \bigm|_{\mathbb C^n\!/\ker(A-\rho_1)}\bigr),\\
 \bar V:=\IM(\bar G''_1+\cdots+\bar G''_p+\rho_1-\rho_2)
 &=\IM(A-\rho_2)/(A-\rho_2)\ker(A-\rho_1)\\
 &=\begin{cases}
    \IM(A-\rho_2)/\ker(A-\rho_1)&(\rho_1\ne\rho_2)\\
    \IM(A-\rho_2)&(\rho_1=\rho_2)
  \end{cases}
\end{align*}
and $(\bar G''_1,\dots,\bar G''_p,\bar G''_1+\cdots+\bar G''_p+\rho_1-\rho_2)$
satisfies the conditions corresponding to \eqref{eq:star} and 
\eqref{eq:starstar}.

Moreover we remark that the last claim on the conditions 
corresponding to \eqref{eq:star} and \eqref{eq:starstar}
doesn't necessarily imply that $(\bar G''_1,\dots,\bar G''_p)$ satisfies
the conditions.

Applying Lemma~\ref{lem:mcbyIM} to $mc_{\rho_1}(\bar G''_1,\dots,\bar G''_p,-\bar G''_1-\cdots-\bar G''_p+\rho_2-\rho_1)$, we have
\begin{align*}
&mc_{\rho_1}
(\bar G''_1,\dots,\bar G''_p,-\bar G''_1-\cdots-\bar G''_p+\rho_2-\rho_1)
 \sim (A'_1,\dots,A'_p,A'_{p+1}),\\
 A'_j &= \bordermatrix{
         &    &     &        &&& \underset{\smile}{p+j} \cr
      & \cr
    j\,{\text{\tiny$)$}} &A_{j1}&\cdots&A_{jp}&A_{j1}&\cdots& A_{jj}-\rho_1
       & \cdots & A_{jp}\cr
      & \cr
  }\in M(2n,\mathbb C),\\
A'_{p+1}&=\begin{pmatrix}
               O_n & O_n\\
               -A+\rho_2 & -A+\rho_1+\rho_2
              \end{pmatrix}\in M(2n,\mathbb C)
\end{align*}
for $j=1,\dots,p$. Here $A'_j$ and $A'_{p+1}$ are endomorphisms of the
linear space
\begin{align}
  U:=\bigl\{\begin{pmatrix}
      u \\ v
     \end{pmatrix}\,;\,
      u\in\mathbb C^n,\ v\in\bar V
     \bigr\}.
\end{align}
Since 
\begin{align}
  A'_1+\cdots+A'_{p+1}=
   \begin{pmatrix}
      A & A-\rho_1\\
      -A+\rho_2 & -A+\rho_2+\rho_1
   \end{pmatrix}
\end{align}
and
\begin{multline}
  \begin{pmatrix}
   I_n \\ & A-\rho_1
  \end{pmatrix}
  \begin{pmatrix}
      A & A-\rho_1\\
      -A+\rho_2 & -A+\rho_1+\rho_2
   \end{pmatrix}\\
 =\begin{pmatrix}
   A &I_n \\
   -(A-\rho_1)(A-\rho_2)& -A+\rho_1+\rho_2
  \end{pmatrix}
  \begin{pmatrix}
   I_n \\ & A-\rho_1
  \end{pmatrix}
\end{multline}
and $A-\rho_1: \bar V \overset{\sim}\to \IM(A-\rho_1)(A-\rho_2)$,
we have
\[
mc_{\rho_1}
(\bar G''_1,\dots,\bar G''_p,-\bar G''_1-\cdots-\bar G''_p+\rho_2-\rho_1)
 \sim (\hat A_1,\dots,\hat A_p,\hat A_{p+1})
\]
with
\begin{align}
 \hat A&:=\begin{pmatrix}
   A &I_n \\
   -(A-\rho_1)(A-\rho_2)& -A+\rho_1+\rho_2
  \end{pmatrix}\in\End(\bar V),\label{eq:barA}\\
 \hat V&:=\mathbb C^n\oplus\IM(A-\rho_1)(A-\rho_2)\label{eq:barV}\\
         &\ =\bigl\{\begin{pmatrix}
           u \\ v
          \end{pmatrix}\,;\, u\in\mathbb C^n,\ 
          v \in \IM(A-\rho_1)(A-\rho_2)
          \bigr\}\subset\mathbb C^{2n},\notag\\
 \hat A_j&:=\diag(O_{n_1+\cdots+n_{j-1}},I_{n_j},O_{n_{j+1}+\cdots+n_p+n})
  \hat A\qquad(j=1,\dots,p),\label{eq:barAj}\\
 \hat A_{p+1}&:=\diag(O_n,I_n)\hat A.
\end{align}
Thus we have the following theorem.
\begin{thm}[Extending operation]\label{thm:ext}
Suppose that the Fuchsian system \eqref{eq:Okubo} of ONF
satisfies \eqref{eq:rank_n} and \eqref{eq:Ostar}.
Then for any complex numbers $\rho_1$ and $\rho_2$ with $\rho_1\rho_2\ne0$,
$(\hat A_1,\dots,\hat A_{p+1}):=
mc_{\rho_1}\circ M^{p+1}_{(0,\dots,0,\rho_2-\rho_1)}\circ T_{(p+1,\infty)}
\circ mc_{-\rho_1}
(A_1,\dots,A_p)$ defines a 
Fuchsian system
\begin{equation}
  \bigl(xI_{\hat n}-\hat T\bigr)\frac{du}{dx} = \hat A u
  \label{eq:ExtO}
\end{equation}
of ONF satisfying \eqref{eq:rank_n} and \eqref{eq:Ostar}.
Here $\hat T=\diag(t_1I_{n_1},\dots,t_pI_{n_p},t_{p+1}I_{n_{p+1}})
\in\End(\hat V)$, $\hat V\simeq\mathbb C^{\hat n}$ and 
$\hat A\in\End(\hat V)$ 
are defined by \eqref{eq:barA} and \eqref{eq:barV} and
\begin{equation}
 \hat n=\dim \hat V=n+n_{p+1},\quad n_{p+1}=\dim\IM(A-\rho_1)(A-\rho_2).
\end{equation}
Moreover \eqref{eq:ExtO} is linearly irreducible if and only if
\eqref{eq:Okubo} is linearly irreducible.

Let
\begin{equation}\label{eq:ORS}
 \begin{Bmatrix}
  x=\infty & x=t_1      & \cdots &x=t_p\\
  [-\mu_1]_{(m_1)}&[0]_{(n-n_1)} & \cdots & [0]_{(n-n_p)}\\
  [-\mu_2]_{(m_2)}&[\lambda_{1,1}]_{(\ell_{1,1})}&\cdots
    &[\lambda_{p,1}]_{(\ell_{p,1})}\\
  \vdots   &\vdots&&\vdots\\
  [-\mu_q]_{(m_q)}&[\lambda_{1,r_1}]_{(\ell_{1,r_1})}&\cdots
    &[\lambda_{p,r_q}]_{(\ell_{p,r_q})}
 \end{Bmatrix}
\end{equation}
be the Riemann scheme of the system \eqref{eq:Okubo} of ONF,
which is compatible with the notation in \eqref{eq:Aii} and \eqref{eq:A}
etc.\ when $A_{ii}$ and $A$ are diagonalizable.
We may assume
\begin{equation}
 \begin{cases}
  \rho_1=\mu_1\text{ \ and \ }\rho_2=\mu_2,\\
  \mu_\nu=\rho_1 \ \Rightarrow \ m_\nu\le m_1,\\
  \mu_\nu=\rho_2 \text{ and }\nu>1\ \Rightarrow \ m_\nu\le m_2.
 \end{cases}
\end{equation}
Here $m_1$ and $m_2$ may be $0$.
Then the Riemann scheme of the system \eqref{eq:ExtO}
equals
\begin{equation}\label{eq:ORSE}
 \begin{Bmatrix}
  x=\infty & x=t_1 & \cdots\! & x=t_p & x=t_{p+1}\\
  [-\mu_1]_{(n-m_2)}&[0]_{(\hat n-n_1)} & \cdots\!& [0]_{(\hat n-n_p)}&[0]_{(n)}\\
  [-\mu_2]_{(n-m_1)}&[\lambda_{1,1}]_{(\ell_{1,1})}&\cdots\!
    &[\lambda_{p,1}]_{(\ell_{p,1})}&[\mu_1+\mu_2-\mu_3]_{(m_3)}\\
   &\vdots   &&\vdots&\vdots\\
    &[\lambda_{1,r_1}]_{(\ell_{1,r_1})}\!\!&\cdots\!
    &[\lambda_{p,r_q}]_{(\ell_{p,r_q})}\!&[\mu_1+\mu_2-\mu_q]_{(m_q)}
 \end{Bmatrix}
\end{equation}
with $\hat n=2n-m_1-m_2$.
\end{thm}
\begin{rem}
{\rm i)\ } Suppose that the system \eqref{eq:Okubo} satisfies 
\eqref{eq:rank_n} and \eqref{eq:Ostar}.
Then 
\begin{align}
  q&\ge2,\\
  \mu_\nu&\ne 0&&(\nu=1,\dots,q),\label{eq:OC0}\\
  \ell_{j,\nu}&\le n-n_j&&(\nu=1,\dots,r_j,\ j=1,\dots,p),\label{eq:OC1}\\
  m_\nu&\le\min\{n_1,\dots,n_p\}&&(\nu=1,\dots,q)\label{eq:OC2}
\end{align}
under the notation in the Theorem~\ref{thm:ext}.
For example the condition 
$\ker(A_j-\lambda_{j,\nu})\cap\bigcap_{\nu\ne j}\ker A_\nu
=\{0\}$ with $\dim\ker A_\nu=n_\nu$ assures \eqref{eq:OC1}.

{\rm ii)\ }
Yokoyama \cite{Yo} defines the extending operation for generic parameters
$\lambda_{j,\nu}$, $\mu_\nu$, $\rho_1$ and $\rho_2$.
It is assumed there that 
$A_{ii}$, $A$, $\hat A_{ii}$ and $\hat A$ are diagonalizable,
$\rank A_{ii}=n_i$, $\rho_1\ne\rho_2$ etc.
In this note we don't assume these conditions.
\end{rem}

{\rm iii)\ }
Applying the extending operation to the equation 
$(x-t_1)\frac{du}{dx}=\lambda u$ with the Riemann scheme 
$\begin{Bmatrix}x=\infty &x=t_1\\-\lambda & \lambda\end{Bmatrix}$, we have 
a Gauss hypergeometric system with the Riemann scheme
$\begin{Bmatrix}x=\infty& x=t_1     & x=t_2\\
                -\rho_1 & 0       & 0\\
                -\rho_2 & \lambda & \rho_1+\rho_2-\lambda
 \end{Bmatrix}$,
which is linearly irreducible.
Here $\lambda$, $\rho_1$ and $\rho_2$ are any complex numbers
satisfying $\rho_1\rho_2\lambda(\rho_1-\lambda)(\rho_2-\lambda)\ne0$.

Theorem~\ref{thm:ext} follows from Theorem~\ref{thm:mid} and 
the argument just before Theorem~\ref{thm:ext}.
We will examine the Riemann scheme of \eqref{eq:ExtO}.
In fact Theorem~\ref{thm:mid} proves that the operation
$M^{p+1}_{(0,\dots,0,\rho_2-\rho_1)}\circ T_{(p+1,\infty)}\circ mc_{-\rho_1}$ 
transforms \eqref{eq:ORS} to
\[
 \begin{Bmatrix}
  x=\infty & x=t_1 & \cdots\! & x=t_p & x=t_{p+1}\\
  [\rho_1-\rho_2]_{(n-m_1)}\!\!\!\!&[0]_{(n-n_1-m_1)} & \cdots\! 
    & [0]_{(n-n_p-m_1)}&[\rho_2-\mu_2]_{(m_2)}\!\\
    &[\lambda_{1,1}-\rho_1]_{(\ell_{1,1})}&\cdots\!
    &[\lambda_{p,1}-\rho_1]_{(\ell_{p,1})}&[\rho_2-\mu_3]_{(m_3)}\!\\
   &\vdots   &&\vdots&\vdots\\
    &[\lambda_{1,r_1}-\rho_1]_{(\ell_{1,r_1})}&\cdots\!
    &[\lambda_{p,r_q}-\rho_1]_{(\ell_{p,r_q})}&[\rho_2-\mu_q]_{(m_q)}
 \end{Bmatrix}
\]
and then the farther operation $mc_{\rho_1}$ to this
gives \eqref{eq:ORSE} because $\rho_2-\mu_2=0$ and 
$\rho_1-\rho_2\ne\rho_1$.
\section{Yokoyama's restricting operation}
Yokoyama's restriction is the inverse of his extension and we have the 
following theorem.
\begin{thm}[Restricting operation]\label{thm:rest}
Let \eqref{eq:Okubo} be a linearly irreducible Fuchsian system of ONF.
Under the notation in Theorem~\ref{thm:ext} we assume $q=2$ and
\begin{equation}\label{eq:CR}
 \mu_1+\mu_2\ne\lambda_{p,\nu}\qquad(\nu=1,\dots,r_p).
\end{equation}
Then 
$mc_{\mu_1}\circ T_{(p,\infty)}
\circ M^p_{(0,\dots,0,\mu_1-\mu_2)}\circ mc_{-\mu_1}(A_1,\dots,A_p)$ 
defines a linearly irreducible Fuchsian system
\begin{equation}
  \bigl(xI_{\check n}-\check T\bigr)\frac{du}{dx} = \check A u\label{eq:Red0}
\end{equation}
of ONF, 
whose Riemann scheme is
\begin{equation}\label{eq:ORR}
 \begin{Bmatrix}
  x=\infty & x=t_1      & \cdots &x=t_{p-1}\\
  [-\mu_1]_{(m_1-n_p)}&[0]_{(\check n-n_1)} & \cdots & [0]_{(\check n-n_{p-1})}\\
  [-\mu_2]_{(m_2-n_p)}&[\lambda_{1,1}]_{(\ell_{1,1})}&\cdots
    &[\lambda_{p-1,1}]_{(\ell_{p-1,1})}\\[-2pt]
   [\lambda_{p,1}-\mu_1-\mu_2]_{(\ell_{p,1})}&\vdots&&\vdots\\[-6.5pt]
   \vdots&\vdots&&\vdots\\
   [\lambda_{p,r_p}-\mu_1-\mu_2]_{(\ell_{p,r_p})}&[\lambda_{1,r_1}]_{(\ell_{1,r_1})}&\cdots
    &[\lambda_{p-1,r_{p-1}}]_{(\ell_{p-1,r_{p-1}})}
 \end{Bmatrix}.
\end{equation}
Here the rank of the resulting system equals 
$\check n=n-n_p=n_1+\cdots+n_{p-1}$ and
\begin{equation}\label{eq:OR2}
 \check T=\begin{pmatrix}
         t_1I_{n_1}\\
         &\ddots\\
         &&t_{p-1}I_{n_{p-1}}
        \end{pmatrix},\quad
 \check A=\begin{pmatrix}
         A_{11}&\cdots&A_{1,p-1}\\
         \vdots&\vdots&\vdots\\
         A_{p-1,1}&\cdots&A_{p-1,p-1}
        \end{pmatrix}.
\end{equation}
\end{thm}
\begin{proof}
Suppose $q=2$.
The operation $M^p_{(0,\dots,0,\mu_1-\mu_2)}\circ mc_{-\mu_1}$
transforms \eqref{eq:ORS} to
\[
 \begin{Bmatrix}
  x=\infty & x=t_1 & \cdots\! & x=t_{p-1} & x=t_p\\
  [0]_{(m_2)}\!\!\!&[0]_{(n-n_1-m_1)} & \cdots\! 
    & [0]_{(n-n_{p-1}-m_1)}&[\mu_1-\mu_2]_{(n-n_p-m_1)}\\
    &[\lambda_{1,1}-\mu_1]_{(\ell_{1,1})}&\cdots\!
    &[\lambda_{p-1,1}-\mu_1]_{(\ell_{p-1,1})}
      &[\lambda_{p,1}-\mu_2]_{(\ell_{p,1})}\\
    &\vdots&&\vdots&\vdots\\
    &[\lambda_{1,r_1}-\mu_1]_{(\ell_{1,r_1})}\!&\cdots\!
    &[\lambda_{p-1,r_q}-\mu_1]_{(\ell_{p-1,r_{p-1}})}\!\!
      &[\lambda_{p,r_q}-\mu_2]_{(\ell_{p,r_p})}
 \end{Bmatrix}
\]
and the farther application $mc_{\mu_1}\circ T_{(p,\infty)}$ to the above
gives \eqref{eq:ORR} because
$\mu_1\ne \mu_1-\mu_2$ and $\mu_1\ne\lambda_{p,\nu}-\mu_2$
for $\nu=1,\dots,r_p$, which corresponds to a system of ONF
as is claimed in Lemma~\ref{lem:mcbyIM}.
Here we note that the rank of the resulting system equals
\begin{align*}
&m_2-\bigl((n-n_1-m_1)+\cdots+(n-n_{p-1}-m_1)+0-(p-2)m_2\bigr)\\
&\quad=n_1+\cdots+n_{p-1}-(p-1)n+(p-1)(m_1+m_2)\\
&\quad=n-n_p.
\end{align*}
Since the restricting operation defined in the theorem is the inverse of the 
extending operation in Theorem~\ref{thm:ext}, we have \eqref{eq:OR2}.
\end{proof}
\begin{rem}\label{rem:generic}
Suppose \eqref{eq:CR} is not valid.
If we apply $E_\tau$ with generic $\tau\in\mathbb C$ to the original system
of ONF preceding to the restriction, the resulting restriction satisfies 
\eqref{eq:CR}.
Note that $mc_\tau$ corresponding to the transformations of
$A$, $\lambda_{j,\nu}$ and $\mu_k$ to
$A+\tau$, $\lambda_{j,\nu}+\tau$ and $\mu_k+\tau$, respectively
(cf.~Corollary~\ref{cor:Euler}).
\end{rem}
\begin{rem} {\rm i) }
The extension and restriction give transformations
between linearly irreducible systems of ONF.
These operations do not change their indices of rigidity.
%

{\rm ii) }
The system \eqref{eq:Okubo} is called \textsl{strongly reducible}
by \cite{Yo} if there exists a non-trivial proper subspace of 
$\mathbb C^n$ which is invariant under $T$ and $A$. 
It is shown there that if the system is not strongly 
reducible, this property is kept by these operations.
\end{rem}

\section{Equivalence of algorithms}
In this section the system \eqref{eq:Schl} of SCF defined by
$\mathbf A=(A_1,\dots,A_p)\in M(n,\mathbb C)^p$ is identified
with the system defined by $\mathbf B\in M(n,\mathbb C)^p$ if 
$\mathbf A\sim \mathbf B$ and then the system is ONF if 
a representative of $\mathbf A$ has the form \eqref{eq:OkSh}.
\begin{prop}\label{prop:COk}
Let $\mathbf A=(A_1,\dots,A_p)\in M(n,\mathbb C)^p$
with \eqref{eq:star} and \eqref{eq:starstar}.
Then $mc_\lambda(\mathbf A)$ is of ONF if and only if
$\lambda$ is not the eigenvalue of $A_0:=-A_1-\cdots-A_p$.
In this case the corresponding system of ONF is given by
\eqref{eq:2Okubo}.
\end{prop}
\begin{proof}
Putting $d=\dim\ker A_1+\cdots+\dim\ker A_p+\dim\ker(A_0-\lambda)
-(p-1)n$, the rank of the system defined by 
$mc_\lambda(\mathbf A)$ equals $n-d$.
Lemma~\ref{lem:Okubo} implies that
$mc_\lambda(\mathbf A)$ is of ONF if and only if
$\sum_{j=1}^p\bigl(n-\dim\ker A_j\bigr)=n-d$, which means
$\dim\ker(A_0-\lambda)=0$.
\end{proof}
\begin{defn}
We denote by $E^p_{\rho_1,\rho_2}$ the extending operation of the system of ONF
given in Theorem~\ref{thm:ext} and by $R^p$ the restricting operation
given in Theorem~\ref{thm:rest}.
Then the restricting operation $R^p_j$ is defined by 
$R^p\circ T_{(j,p)}$ for $j=1,\dots,p$.
Here $(j,p)$ is the transposition of indices $j$ and $p$
(cf.~Definition~\ref{def:T}).
Note that the restricting operation is defined only when $q=2$.
\end{defn}

We have proved that the extension and the restriction
of the system of ONF is realized by suitable 
combinations of additions, middle convolutions and the automorphism
of $\mathbb P^1(\mathbb C)$ written by $T_{(p+1,\infty)}$ and $T_\sigma$
(cf.~Definition~\ref{def:T}).

In fact, we have the following equalities for operations to linearly 
irreducible systems \eqref{eq:Okubo} of ONF.
\begin{gather}
 E^p_{\rho_1,\rho_2} = mc_{\rho_1}\circ M^{p+1}_{(0,\dots,0,\rho_2-\rho_1)}
 \circ T_{(p+1,\infty)}\circ mc_{-\rho_1},\label{eq:Emc}\\
 R^p = mc_{\mu_1}\circ T_{(p,\infty)}\circ M^p_{(0,\ldots,0,
   {\mu_1-\mu_2})}
   \circ mc_{-\mu_1},\label{eq:Rmc}\\
  R^{p+1}\circ E^p_{\rho_1,\rho_2}= \id.
\end{gather}
Here $\mu_1$ and $\mu_2\in\mathbb C$ are determined by
\begin{align}
(A-\mu_1)(A-\mu_2)=0.
\end{align}
\begin{lem}\label{lem:ERM}
We have the following relations for $j=1,\dots,p$.
\begin{align}
 &R^{p+1}_j \circ E_{\epsilon} \circ E^p_{\rho_1,\rho_2}
 =mc_{\rho_1+\epsilon}
  \circ M^p_{(0,\ldots,0,
   \underset{\overset{\smallfrown}{j}}{\rho_2-\rho_1},
   0,\ldots,0)}
  \circ T_{(j,\infty)}\circ mc_{-\rho_1},\label{eq:RE}\\
 &\ord R^{p+1}_j \circ E_{\epsilon} \circ E^p_{\rho_1,\rho_2}(\mathbf A)
  = \ord \mathbf A+\dim \IM(A-\rho_1)(A-\rho_2)-\dim\IM A_j,
    \label{eq:ordER}\\
 & R^{p+1}_j\circ E_{\rho_1+\epsilon,
   \rho_1+\rho_2+\rho_3+\epsilon}\circ 
   R^{p+1}_j\circ E_\epsilon\circ E^p_{\rho_1,\rho_2}\label{eq:RERE}\\
 &\quad
 =mc_{\rho_1+\epsilon}\circ
   M^p_{(0,\ldots,0,
    \underset{\overset{\smallfrown}{j}}{\rho_1+\rho_3},
    0,\ldots,0)}
 \circ mc_{-\rho_1}.\notag
\end{align}
Here $\rho_1$ and $\rho_2$ are any non-zero complex numbers
and $\epsilon$ is a generic complex number 
(cf.~Remark~\ref{rem:generic}) and $\ord\mathbf A$ denotes the rank of the
corresponding system \eqref{eq:Schl} of SCF.
\end{lem}
\begin{proof}
We may assume $j=1$. It follows from \eqref{eq:Emc} and \eqref{eq:Rmc}
that
\begin{align*}
 &R^{p+1}_1 \circ E_{\epsilon} \circ E^p_{\rho_1,\rho_2}\\
 &= mc_{\rho_1+\epsilon}\circ T_{(p+1,\infty)}
  \circ M^{p+1}_{(0,\ldots,0,\rho_1-\rho_2)}\circ
  mc_{-\rho_1-\epsilon}\circ T_{(1,p+1)}\circ
  mc_{\epsilon}\\
 &\quad \circ mc_{\rho_1}\circ M^{p+1}_{(0,\ldots,0,\rho_2-\rho_1)}
 \circ T_{(p+1,\infty)} \circ mc_{-\rho_1}\\
 &=mc_{\rho_1+\epsilon}\circ T_{(p+1,\infty)}
  \circ M^{p+1}_{(0,\ldots,0,\rho_1-\rho_2)}\circ T_{(1,p+1)}
  \circ M^{p+1}_{(0,\ldots,0,\rho_2-\rho_1)}
 \circ T_{(p+1,\infty)}\\
 &\quad \circ mc_{-\rho_1}\\
 &=mc_{\rho_1+\epsilon}\circ T_{(p+1,\infty)}
  \circ M^{p+1}_{(\rho_2-\rho_1,0,\ldots,0,\rho_1-\rho_2)}
 \circ T_{(1,p+1)}\circ T_{(p+1,\infty)} \circ mc_{-\rho_1}\\
 &=mc_{\rho_1+\epsilon}
  \circ M^p_{(\rho_2-\rho_1,0,\ldots,0)}
  \circ T_{(1,\infty)}\circ mc_{-\rho_1}
\intertext{and therefore}
 &R^{p+1}_1\circ E^p_{\rho_1+\epsilon,
 \rho_1+\rho_2+\rho_3+\epsilon}\circ 
 R^{p+1}_1\circ E_\epsilon\circ E^p_{\rho_1,\rho_2}\\
 &=mc_{\rho_1+\epsilon}\circ
   M^p_{(\rho_2+\rho_3,0,\ldots,0)}\circ T_{(1,\infty)}\circ
   mc_{-\rho_1-\epsilon}
 \\&\quad\circ
 mc_{\rho_1+\epsilon}
 \circ M^p_{(-\rho_1+\rho_2,0,\ldots,0)}
 \circ T_{(1,\infty)} \circ mc_{-\rho_1}\\
 &=mc_{\rho_1+\epsilon}\circ
   M^p_{(\rho_1+\rho_3,0,\ldots,0)}\circ mc_{-\rho_1}.
\end{align*}
The equality \eqref{eq:ordER} follows from Theorem~\ref{thm:ext} and
Theorem~\ref{thm:rest}.
\end{proof}
We show Riemann schemes related to \eqref{eq:RERE}.

\begin{rem}\label{rem:RERE}
By the extension $E^p_{-\lambda_{0,1},-\lambda_{0,2}}$ we have
{\small
\begin{multline*}
\begin{Bmatrix}
  x=\infty & x=t_1&\cdots & x=t_p\\
  [\lambda_{0,1}]_{(m_{0,1})} & [0]_{(m_{1,1})}& \cdots&[0]_{(m_{p,1})}\\
  [\lambda_{0,2}]_{(m_{0,2})} & [\lambda_{1,2}]_{(m_{1,2})} & \cdots
   &[\lambda_{p,2}]_{(m_{p,2})}\\
  \vdots & \vdots & & \vdots
\end{Bmatrix}\mapsto
\\
\begin{Bmatrix}
  x=\infty & x=t_1 & \cdots & x=t_{p+1}\\
  [\lambda_{0,1}]_{(n-m_{0,2})}& [0]_{(m_{1,1}+n-m_{0,1}-m_{0,2})}
  & \cdots & [0]_{(n)}\\
  [\lambda_{0,2}]_{(n-m_{0,1})}& [\lambda_{1,2}]_{(m_{1,2})} &\cdots&
  [\lambda_{0,3}-\lambda_{0,1}-\lambda_{0,2}]_{(m_{0,3})} \\
  &\vdots & & \vdots
\end{Bmatrix}.
\end{multline*}}

\noindent
Here $n=m_{j,1}+m_{j,2}+\cdots$ and $m_{1,1}+\cdots+m_{p,1}=(p-1)n$.
By applying the restriction $R^{p+1}_1\circ E_\epsilon$ to this result 
we have
{\small \begin{align*}
\begin{Bmatrix}
x=\infty & x=t_1 & x=t_2 & \cdots\\
 [\lambda_{0,1}-\epsilon]_{(m_{1,1}-m_{0,2})}\!\!
  & [0]_{(m_{1,1})}
   & [0]_{(m_{2,1}-m_{0,1}-m_{0,2}+m_{1,1})}\!\!\!\! & \cdots\\
 [\lambda_{0,2}-\epsilon]_{(m_{1,1}-m_{0,1})}\!\!
  & [\lambda_{0,3}-\lambda_{0,1}-\lambda_{0,2}+\epsilon]_{(m_{0,3})}
   \!\!\!\!\!\!
   &[\lambda_{2,2}+\epsilon]_{m_{2,2}}&\cdots\\[-4pt]
 [\lambda_{1,2}+\lambda_{0,1}+\lambda_{0,2}-\epsilon]_{(m_{1,2})}
  \!\!\!\!\!\!& \vdots& \vdots & \cdots
  \\
[\lambda_{1,3}+\lambda_{0,1}+\lambda_{0,2}-\epsilon]_{(m_{1,3)}}
  \!\!\!\!\!\!\\
\vdots
\end{Bmatrix}
\end{align*}}

\noindent
whose rank equals $n-(m_{0,1}+m_{0,2}-m_{1,1})$.
By applying the extending operation 
$E^p_{-\lambda_{0,1}+\epsilon,-\lambda_{1,2}-\lambda_{0,1}-\lambda_{0,2}+\epsilon}$ to what we obtained we have
{\small\begin{align*}
&\left\{
\begin{matrix}
x=\infty & x=t_1\\
[\lambda_{0,1}-\epsilon]_{(n-m_{0,1}-m_{0,2}+m_{1,1}-m_{1,2})}
 & [0]_{(n-m_{0,1}+m_{1,1}-m_{1,2})}
\\
[\lambda_{1,2}+\lambda_{0,1}+\lambda_{0,2}-\epsilon]_{(n-m_{0,1})}
 &[\lambda_{0,3}-\lambda_{0,1}-\lambda_{0,2}+\epsilon]_{(m_{0,3})}
\\
  &[\lambda_{0,4}-\lambda_{0,1}-\lambda_{0,2}+\epsilon]_{(m_{0,4})}
\\
  &\vdots
\end{matrix}\right.
\\
&\qquad\qquad\left.
\begin{matrix}
x=t_2 & \cdots & x=t_{p+1}
\\
   [0]_{(n-2m_{0,1}-m_{0,2}+m_{1,1}-m_{1,2}+m_{2,1})}
   &\cdots 
    & [0]_{(n-m_{0,1}-m_{0,2}+m_{1,1})}
\\
   [\lambda_{2,2}+\epsilon]_{(m_{2,2})}
   &\cdots 
    &[-2\lambda_{0,1}-\lambda_{1,2}+\epsilon]_{(m_{1,1}-m_{0,1})}
\\[-4pt]
     \vdots
     & 
    \ldots
     &[\lambda_{1,3}-\lambda_{1,2}-\lambda_{0,1}+\epsilon]_{(m_{1,3})}
\\
 & & \vdots
\end{matrix}
\right\}
\end{align*}}

\noindent
and by applying the restriction $R^{p+1}_1$ to this result we finally have
{\small\begin{align*}
\begin{Bmatrix}
x=\infty & x=t_1 & x=t_2 & \cdots \\
[\lambda_{0,1}-\epsilon]_{(m_{0,1}-d)}\!\!\!\!\!
 & [0]_{(m_{1,1})}
  &[0]_{(m_{2,1}-d)}
   &\cdots\\
[\lambda_{0,2}+\lambda_{0,1}+\lambda_{1,2}-\epsilon]_{(m_{0,2})}\!\!
 &[-2\lambda_{0,1}-\lambda_{1,2}+\epsilon]_{(m_{1,1}-m_{0,1})}\!\!\!\!
  &[\lambda_{2,2}+\epsilon]_{(m_{2,2})}\!\!
   &\cdots\\[-4pt]
[\lambda_{0,3}+\lambda_{0,1}+\lambda_{1,2}-\epsilon]_{(m_{0,3})}
 &[\lambda_{1,3}-\lambda_{1,2}-\lambda_{0,1}+\epsilon]_{(m_{1,3})}
  & \vdots
   & \cdots\\
[\lambda_{0,4}+\lambda_{0,1}+\lambda_{1,2}-\epsilon]_{(m_{0,4})}
 &[\lambda_{1,4}-\lambda_{1,2}-\lambda_{0,1}+\epsilon]_{(m_{1,4})}\\
 \vdots&\vdots
\end{Bmatrix}
\end{align*}}

\noindent
with $d=m_{0,1}-m_{1,1}+m_{1,2}$.
\end{rem}
\begin{thm}\label{thm:eq}
Suppose $\mathbf A=(A_1,\dots,A_p)\in M(n,\mathbb C)^p$ 
is irreducible 
and suppose $\mathbf B=(B_1,\dots,B_p)\in M(n,\mathbb C)^p$ 
is obtained from $\mathbf A$ by a finite iteration of additions,
middle convolutions and operations $T_{(p,\infty)}$ and $T_\sigma$ 
in Definition~\ref{def:T}.

Let $\alpha$ and $\beta$ be generic complex numbers
so that $mc_\alpha(\mathbf A)$ and $mc_\beta(\mathbf B)$ are of ONF.
Then $mc_\beta(\mathbf B)$ can be obtained from $mc_\alpha(\mathbf A)$ 
by a finite iteration of the suitable operations
$R^{p+1}_j\circ E_\epsilon\circ E^p_{\rho_1,\rho_2}$, namely,
extensions, restrictions and generic Euler transformations.
Here $\alpha=0$ is generic if $\mathbf A$ is of ONF.
\end{thm}
\begin{proof}
The theorem follows from Lemma~\ref{lem:ERM} 
since $R_j^{p+1}\circ E^p_{\rho_1,\rho_1}=T_{(j,\infty)}$,
$T_{(i,j)}=T_{(j,\infty)}\circ T_{(i,\infty)}\circ T_{(j,\infty)}$, 
$M_\mu\circ M_{\mu'}=M_{\mu+\mu'}$,
$mc_\lambda\circ mc_{\lambda'}=mc_{\lambda+\lambda'}$
and $mc_0=id$. 
\end{proof}
\section{Reduction process}
For the system \eqref{eq:Schl} of SCF the \textsl{spectral type} of $\mathbf A=
(A_1,\dots,A_p)$ denoted by $\spt\mathbf A$ is the tuple of 
$p+1$ partitions of $n$
\begin{equation}
 \spt\mathbf A:=\mathbf m = 
   (m_{0,1},\dots,m_{0,n_0};m_{1,1},\dots,m_{1,n_1};\dots;
   m_{p,1},\dots,m_{p,n_p})
\end{equation}
under the notation \eqref{eq:DS}.
This tuple may be expressed by
\begin{equation}
m_{0,1}\cdots m_{0,n_0},m_{1,1}\cdots m_{1,n_1},\cdots,
   m_{p,1}\cdots m_{p,n_p}
\end{equation}
and in this case \eqref{eq:ZL} shows
\begin{equation}
 \idx\mathbf A=\sum_{\substack{1\le\nu\le n_j\\ 0\le j\le p}}
   m_{j,\nu}^2 - (p-1)(\ord\mathbf A)^2.
\end{equation}
We put $n_j=1$ and $m_{j,1}=\ord\mathbf m:=m_{0,1}+\cdots+m_{0,n_0}$ 
if $j>p$.  Moreover we put $m_{j,\nu}=0$ if $j>n_j$.

For $p+1$ non-negative integers $\tau=(\tau_0,\dots,\tau_p)$
we define
\begin{equation}
 d_\tau(\mathbf m) := m_{0,\tau_0}+\cdots+m_{p,\tau p} - (p-1)\ord\mathbf A
\end{equation}
and $\tau(\mathbf m)=\bigl(\tau(\mathbf m)_0,\dots,\tau(\mathbf m)_p\bigr)$
so that
\begin{equation}
 m_{j,\tau(\mathbf m)_j}\ge m_{j,\nu}\qquad(\nu=1,\dots,n_j,\ j=0,\dots,p).
\end{equation}
Moreover we put
\begin{align}
 d_{\max}(\mathbf m):=d_{\tau(\mathbf m)}(\mathbf m).
\end{align}

Suppose $\mathbf A$ is irreducible.
Put
\begin{align}
  mc_{\max}(\mathbf A):=
  mc_{\lambda_{0,\tau(\mathbf m)_0}+\cdots+\lambda_{p,\tau(\mathbf m)_p}}
   \circ M_{(-\lambda_{1,\tau(\mathbf m)_1},\dots,
   -\lambda_{p,\tau(\mathbf m)_p})}(\mathbf A)
\end{align}
under the notation \eqref{eq:DS}.
If $n>1$, then Theorem~\ref{thm:mid} proves
\begin{align}
 &\begin{cases}
  \spt mc_{max}(\mathbf A)=
  \p_{max}(\mathbf m):=(\ldots;m'_{j,1},\ldots,m'_{j,n_j};\ldots)\\
   m'_{j,\nu}= m_{j,\nu}-d_{\max}(\mathbf m)\delta_{\nu,\tau(\mathbf m)_j}
   \qquad(\nu=1,\dots,n_j,\ j=0,\dots,p)
 \end{cases},\\
 &\quad\ord \p_{max}(\mathbf m) = \ord\mathbf m-d_{max}(\mathbf m).
\end{align}

If $\mathbf A$ is rigid,
namely, $\idx\mathbf m=2$, then we have $d_{max}(\mathbf m)>0$
because
\begin{equation*}
 \idx\mathbf m +  \sum_{j=0}^p\sum_{\nu=1}^{n_j}
 (m_{j,\tau_j}-m_{j,\nu})\cdot m_{j,\nu}
=\biggl(\sum_{j=0}^pm_{j,\tau_j}-(p-1)\ord\mathbf m\biggr)
 \cdot\ord\mathbf m
\end{equation*}
and thus we have 
$\ord mc_{max}(\mathbf A)<\ord\mathbf A$.
Hence if the system of SCF is linearly irreducible and rigid,
the system is connected to a rank 1 system by a finite iteration 
of additions and middle convolutions and conversely any linearly irreducible
system of SCF is constructed from a rank 1 system by a finite iteration of 
additions and middle convolutions
(cf.~\cite{Kz}, \cite{Ko}, \cite{DR}, \cite{O2}).

Since any rank 1 system is transformed into ONF by a suitable addition, 
Theorem~\ref{thm:eq} implies the following theorem, which is 
given in \cite[Theorem~4.6]{Yo} when the parameters
$\lambda_{i,\nu}$ and $\mu_j$ are generic.

\begin{thm}\label{thm:Oreg}
Any linearly irreducible rigid system of ONF is connected to a rank $1$ system 
of ONF by a finite iteration of extensions, restrictions and generic Euler 
transformations.
\end{thm}

\begin{rem} {\rm i) }
For a given $\mathbf A\in M(n,\mathbb C)^p$, if there exists
$j$ with  $d_{max}(\spt\mathbf A)>m_{j,\tau(\spt\mathbf A)_j}$, 
$\mathbf A$ is not irreducible. This is a consequence of
Theorem~\ref{thm:ext}.

{\rm ii) }
It follows from Proposition~\ref{prop:COk} that $mc_{\max}(\mathbf A)$ is 
not of ONF for any linearly irreducible system \eqref{eq:Schl} of SCF.

{\rm iii) }
In virtue of Lemma~\ref{lem:ineqforKatz} 
a more explicit construction of the reduction process within
ONF using extensions, restrictions and generic Euler
transformations is obtained as follows. 

Put $\mathbf m=\spt(\mathbf A)$ for a linearly irreducible system
\eqref{eq:Okubo} of ONF.
Assume that $\mathbf m$ satisfies the assumption of 
Lemma~\ref{lem:ineqforKatz} and $\lambda_{j,1}= 0$ for $j=1,\dots,p$.
Then Lemma~\ref{lem:ineqforKatz} assures that we can find $j\ge 1$ 
with
\begin{equation}\label{eq:ineqOk}
 m_{0,1}-m_{j,1}+m_{j,2}>0
\end{equation}
because $d_{max}(\mathbf m)=m_{0,1}$.
Applying the operation \eqref{eq:RERE} with $\rho_1=\lambda_{0,1}$,
$\rho_2=\lambda_{0,2}$ and $\rho_3=\lambda_{1,2}$, it follows from
Remark~\ref{rem:RERE} that
the resulting $\mathbf A'$ satisfies 
\begin{equation}
 \ord\mathbf A'=\ord\mathbf A-m_{0,1}+m_{j,1}-m_{j,2} < \ord\mathbf A.
\end{equation}

{\rm iv)\ }
The existence of $j\ge 1$ satisfying \eqref{eq:ineqOk} is 
given by \cite[Lemma~4.2]{Yo} when the rigidity index of 
the system of ONF equals $2$.
Note that any linearly irreducible rigid system of SCF with rank $> 1$ always
satisfies the assumption of Lemma~\ref{lem:ineqforKatz}.
\end{rem}
\begin{lem}\label{lem:ineqforKatz}
\iffalse
Let $\mathbf m=(\mathbf m_0,\dots,\mathbf m_p)$ 
be a tuple of $p+1$ partitions of 
$\mathbf m_j=(m_{j,1},\dots,m_{j,n_j})$ of $n$
satisfying $m_{j,2}\ge m_{j,3}\ge\cdots\ge m_{j,n_j}$ for $j=0,\dots,p$.
Suppose $d:=m_{0,1}+\cdots+m_{p,1}-(p-1)\ord\mathbf m > 0$,
and put $\mathbf m'=(\ldots;m_{j,1}-d, m_{j,2},\dots,m_{j,n_j};\ldots)$
If $d_{max}(\mathbf m')>0$, then
\else
Let $\mathbf m$ be a spectral type of a linearly irreducible system
\eqref{eq:Schl} of SCF with $\ord\mathbf m>1$.
Put $\mathbf m'=\p_{max}(\mathbf m)$.
We may assume $m_{j,1}\ge m_{j,2}\ge\cdots \ge m_{j,n_j}$.
If $d_{max}(\mathbf m)>0$ and $d_{max}(\mathbf m')>0$, then
\fi
\begin{equation}\label{eq:ineqKz}
\sum_{j=0}^p\max\{0,d_{max}(\mathbf m)
  -(m_{j,1}-m_{j,2})\} > d_{max}(\mathbf m).
\end{equation}
\end{lem}
\begin{proof}
Put $d=d_{max}(\mathbf m)$.
Since $\max\{m'_{j,1},\dots,m'_{j,n_j}\}=
\max\{m_{j,2}, m_{j,1}-d\}$, the assumption implies
\[
  \sum_{j=0}^p \max\{m_{j,2},m_{j,1}-d\} > 
  (p-1)\ord\mathbf m'=(p-1)(n-d).
\]
Hence we have
\begin{align*}
 \sum_{j=0}^p \max\{d-(m_{j,1}-m_{j,2}),0\} &>
 (p-1)(n-d) - \sum_{j=1}^p(m_{j,1}-d)\\
 &=(p-1)(n-d) - (p-1)n + pd = d.
\end{align*}\\[-.97cm]
\end{proof}

A linearly irreducible system \eqref{eq:Schl} of SCF satisfying
$d_{max}(\spt\mathbf A)\le 0$ is called \textsl{basic}, 
which is not rigid and not of ONF.
It is known that the basic systems of SCF
with different spectral types cannot be connected by any iteration
of middle convolutions, additions, $T_{j,\infty}$ and $T_\sigma$.
Moreover there exist a finite number
of basic systems with a fixed index of rigidity and an indivisible
spectral type (cf.~\cite{CB}, \cite[Proposition~8.1]{O2}).
Here $\mathbf m=(\ldots;m_{j,1},\dots,m_{j,n_j};\ldots)$ is 
indivisible if there doesn't exist a non-trivial common divisor
of $\{m_{j,\nu}\,;\,j=0,1,\ldots,\ \nu=1,2,\ldots\}$
and two tuples are identified if a permutation of indices $j$
and permutations of indices $\nu$ within the same $j$ 
transform one of the two into the other.

It is shown by \cite{CB} that the basic systems with a given
index of rigidity correspond to the positive imaginary roots with a fixed
norm in the closure of a negative Weyl chamber of a Kac-Moody root system
with a star-shaped Dynkin diagram (cf.~\cite{Kc}, \cite[\S7]{O2}).
Any linearly irreducible system of SCF which is not rigid
is connected to a basic system by an iteration of $mc_{max}$
and therefore we have the following theorem.
\begin{thm}
By a finite iteration of extensions, restrictions
and generic Euler transformations,
any linearly irreducible system of ONF which is not rigid is connected
to a system of ONF transformed by a middle convolution of a basic 
system of SCF (cf.~Proposition~\ref{prop:COk}).
\end{thm}
We will give some examples.

\begin{exmp}
There exist 4 different spectral types of
basic systems with index of rigidity $0$
(cf.~\cite{Ko2}, \cite[Proposition~8.1]{O2}):

\smallskip
\centerline{\begin{tabular}{|c||c|l||c|l|}\hline
 type         & $\ord$ & basic system & $\ord$ &ONF\\ \hline
 $\tilde D_4$ & 2& 11,11,11,11 & 3&111,21,21,21\\ \hline
 $\tilde E_6$ & 3&111,111,111  & 4&1111,211,211\\ \hline
 $\tilde E_7$ & 4&1111,1111,22 & 5&11111,2111,32\\ \hline
 $\tilde E_8$ & 6&111111,222,33& 7&1111111,322,43\\ \hline
\end{tabular}}

\bigskip
The following is the list of spectral types of basic systems 
with index of rigidity $-2$
(cf.~\cite[Proposition~8.4]{O2}):
\smallskip

\centerline{\begin{tabular}{|r|l||r|l|}\hline 
  $\ord$ &basic system& $\ord$ &ONF \\ \hline
  2& 11,11,11,11,11&4& 211,31,31,31,31\\ \hline
  3& 111,111,21,21 &4& 1111,211,31,31\\ \hline
  4&1111,22,22,31  &5& 11111,32,32,41\\ \hline
  4&1111,1111,211  &5& 11111,2111,311\\ \hline
  4&211,22,22,22   &6& 2211,42,42,42\qquad 222,411,42,42\\ \hline
  5&11111,221,221  &6& 111111,321,321\\ \hline
  5&11111,11111,32 &6& 111111,21111,42\\ \hline
  6&111111,2211,33 &7& 1111111,3211,43\\ \hline
  6&2211,222,222   &8& 22211,422,422\qquad 2222,422,4211\\ \hline
  8&11111111,332,44&9& 111111111,432,54\\ \hline
  8&22211,2222,44  &10&222211,4222,64\quad \ \ 22222,42211,64\\ \hline
 10&22222,3331,55  &12&222222,5331,75\\ \hline
 12&2222211,444,66 &14&22222211,644,86\\ \hline
\end{tabular}}

\medskip
Here we give the spectral types of systems of ONF with 
minimal rank corresponding to a basic system, 
which are not necessarily unique but transformed to each other 
by suitable iterations of extensions, restrictions 
and generic Euler transformations.
\end{exmp}
\begin{defn}
For a $(p+1)$-tuple $\mathbf m=(m_{j,\nu})$ of partitions of $n$ we put
\begin{equation}
  \Oidx\mathbf m:= (p-1)\cdot\ord\mathbf m - \max_{0\le k\le p}
                 \sum_{\substack{0\le j\le p\\j\ne k}}
                 \max\{m_{j,1},m_{j,2},\ldots\}.
\end{equation}
We define that $\mathbf m$ is of \textsl{Okubo type} if 
$\Oidx\mathbf m=0$.
\end{defn}
\begin{rem}
Let $\mathbf m$ be the spectral type of a linearly irreducible
system \eqref{eq:Schl} of SCF.  Then $\Oidx\mathbf m\ge 0$.
Moreover it follows form Lemma~\ref{C:Okubo} and Remark~\ref{rem:C:Okubo} 
that the system is equivalent to a system of ONF after applying a 
suitable addition (and $T_{(p,\infty)}$) if and only if $\Oidx\mathbf m=0$.

If $\mathbf m$ is basic, then $\Oidx\mathbf m> 0$ and 
there exists a system of ONF with the minimal rank
$\ord\mathbf m+\Oidx\mathbf m$ among the systems obtained from 
the original system \eqref{eq:Schl} by a finite iteration
of additions and middle convolutions
(cf.~Proposition~\ref{prop:COk}).
\end{rem}

\end{document}